\documentclass[3p,times]{article}

\usepackage{amsmath}
\usepackage{amsthm}
\usepackage{booktabs}
\newtheorem{prop}{Proposition}

\newtheorem{cor}{Corollary}
\newtheorem{exmp}{Example}
\newtheorem{obs}{Observation}
\newtheorem{rem}{Remark}

\usepackage{amssymb}

\usepackage[figuresright]{rotating}
\usepackage{tcolorbox}

\usepackage[font=footnotesize,labelfont=bf]{caption}
\usepackage[font=footnotesize,labelfont=bf]{subcaption}

\usepackage[english]{babel}
\title{On an Umbral point of view to the Gaussian and Gaussian like functions}

\author{Giuseppe Dattoli, Emanuele Di Palma, Silvia Licciardi\footnote{Corresponding author: silviakant@gmail.com, silvia.licciardi@enea.it, orcid 0000-0003-4564-8866, tel. nr: +39 06-94005421.},\\[1.3ex]
	ENEA - Frascati Research Center, Via Enrico Fermi 45, 00044, \\Frascati, Rome, Italy  
}

\date{\today}
\begin{document}
	
	\maketitle

\abstract{In this note we review the theory of Gaussian functions by exploiting a point of view based on symbolic methods of umbral nature. We introduce quasi-Gaussian functions, which are close to Gaussian distribution but have a longer tail. Their use and their link with hypergeometric function is eventually presented.}


	\section{Introduction}

Umbral Methods ($UM$) \cite{S.Roman,SLicciardi} realize a common environment in which the properties of different functions can be studied, using a specific image realization \cite{L.C.Andrews,Abramovitz}.
The Bessel functions have indeed been studied using the Gaussian as the associated image function \cite{SLicciardi}. The trigonometric functions \cite{FromCircular,Airy,SLicciardi} too have been framed within a context having the Gaussian as a pivot reference function. Within this context a natural transition between circular and Bessel function has been obtained and the spherical Bessel have been shown to be the thread union between the two families of functions.\\

\noindent The paradigm, we have followed, has been that of providing a kind of ``downgrade'' of the level of the complexity of the function itself, by reducing a higher order transcendental function to a lower order.
The simplification occurs through the introduction of a set of operators, with specific algebraic properties, whose role has been clarified and made more rigorous within the well established properties of the Borel transform \cite{Babusci,Prod}.  \\

\noindent The procedure, we have just summarized, has revealed interesting aspects of special function unvealing elements of superposition, hardly achieveble by other means.
In this paper we use the same point of view by exploring the consequences of the ``downgrading" of the Gauss function to the status of a rational function.\\

\noindent The starting point of our study is the reduction of the Gaussian to a Lorentz function and we will see how this opens a view on other families of functions, leading to generaled forms of Gaussian trigonometric like functions.
The paper is organized as outlined below.\\

\noindent In Section 2, we clarify how a Lorentzian can be chosen as the relevant umbral image and explore further associated umbral forms, which naturally leads to the sine and cosine Gaussian functions.\\
\noindent In Section 3 we extend the formalism to the study of higher order trig-like functions and introduce further generalizations.\\
\noindent Finally, Section 4 contains comments, including applications and possible links of the previous conclusions with Lèvy distributions.

\section{Gaussian Functions, Lorentzian and Associated Trigonometric Functions}

We make use of umbral operators to construct images of special and ordinary functions. The underlying formalism revealed quite powerful to deal with computational details (difficult to accomplish with ordinary means) and to disclosing intimate relationships between different forms of special functions.\\

\noindent In this section, we show the consequences deriving from the umbral restyling of the Gauss function in terms of the Lorentz function.

\begin{prop}
	We write the Gauss function as
	\begin{equation} \label{GrindEQ__1_} 
	e^{-x^{2} } =\frac{1}{1+\hat{c}\, x^{2} } \varphi _{0} , 
	\end{equation} 
	where the umbral operator $\hat{c}$ is such that
	
	\begin{equation} \label{GrindEQ__2_} 
	\hat{c}\, ^{\alpha } \varphi _{0} =\frac{1}{\Gamma (\alpha +1)} \; ,
	\end{equation} 
	with $\Gamma (.)$ being the Euler gamma function and $\alpha $ any real or complex number.
\end{prop}
\begin{proof}
	According to the umbral point of view\footnote{See \cite{SLicciardi} for a rigorouss treatment of the Umbral Methods.}
	\begin{equation}
	e^{-x^{2} } =\sum_{r=0}^\infty \dfrac{(-1)^r x^{2r}}{r!}=
	\sum_{r=0}^\infty \dfrac{(-1)^r x^{2r}}{\Gamma(r+1)}=
	\sum_{r=0}^\infty (-\hat{c})^r x^{2r}\varphi_{0}=
	\frac{1}{1+\hat{c}\, x^{2} } \varphi _{0} . 
	\end{equation}
\end{proof}	

\begin{cor}
 A consequence of Eq. \eqref{GrindEQ__1_} is that the primitive of the Gaussian function too can be formally expressed in terms of an elementary function, namely\footnote{We use the notation $\tan^{-1}(x)$ to indicate $\arctan(x)$.}

\begin{equation}
\int e^{-x ^{2} }  dx =\hat{c}^{-\frac{1}{2} } tg^{-1} \left(\hat{c}^{\frac{1}{2} } x\right) .
\end{equation}  
\end{cor}
           
\begin{figure}[h]
	\centering
	\includegraphics[width=0.4\linewidth]{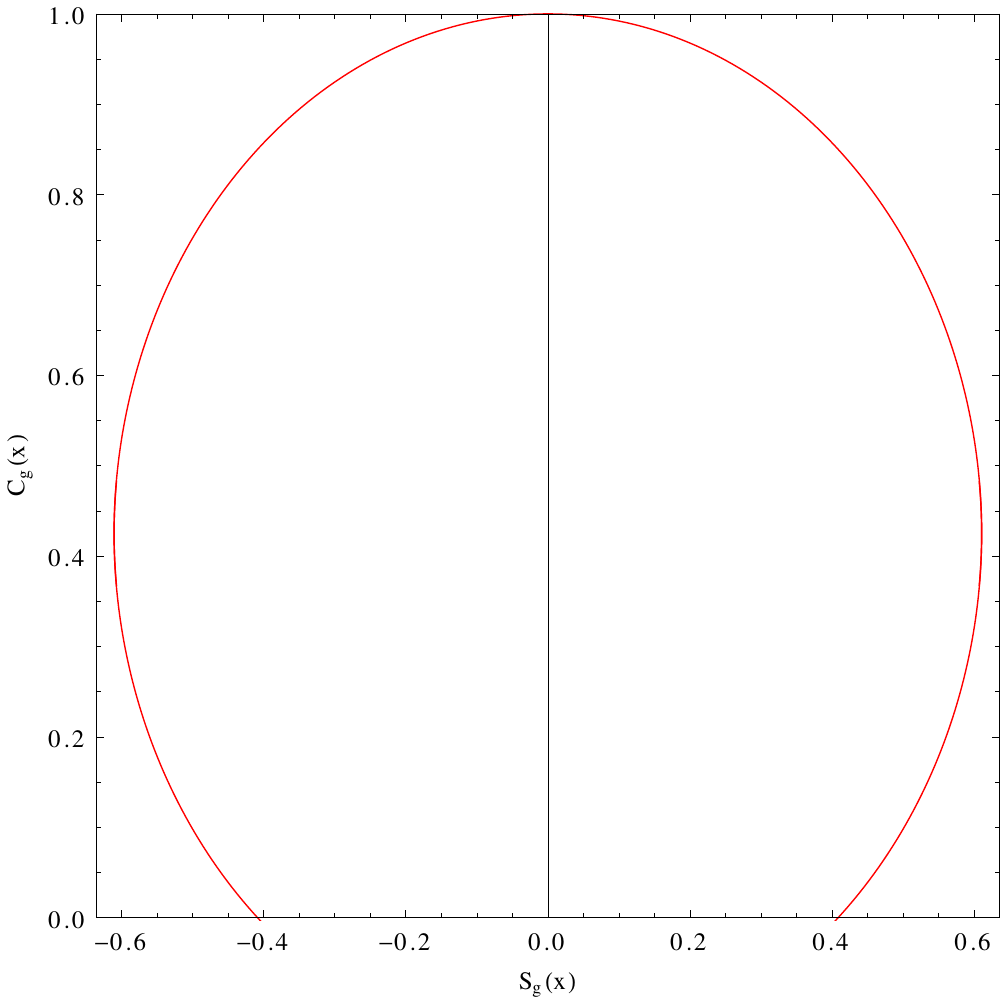}
	\caption{Gaussian trigonometric "circumference": egg-shaped curve. $C_{g} (x)$ vs $S_{g} (x)$. }\label{fig1} 
\end{figure}

\begin{figure}[h]
	\centering
	\includegraphics[width=0.5\linewidth]{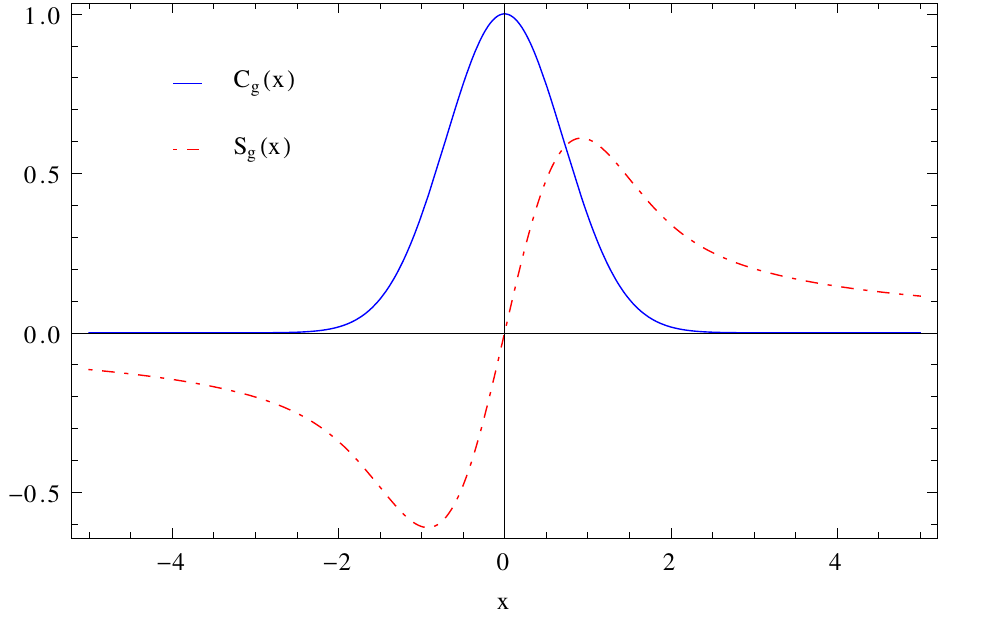}
	\caption{$C_{g} (x)$ and $S_{g} (x)$ vs. $x$.}\label{fig2} 
\end{figure}

 \begin{prop}          
 It is spossible cast Eq. \eqref{GrindEQ__1_} in the form
\begin{equation} \label{GrindEQ__4_} 
e^{-x^{2} } :=C_{g} (x)=\frac{1}{2} \left(\frac{1}{1-i\, \hat{c}^{\frac{1}{2} } x} +\frac{1}{1+i\, \hat{c}^{\frac{1}{2} } x} \right)\, \varphi _{0}  
\end{equation} 
and define the associated function
\begin{equation}\label{GrindEQ__5_} 
S_{g} (x):=\frac{1}{2\, i} \left[\frac{1}{1-i\, \hat{c}^{\frac{1}{2} } x} -\frac{1}{1+i\, \hat{c}^{\frac{1}{2} } x} \right]\, \varphi _{0} = \frac{2}{\sqrt{\pi } } \sum _{r=0}^{\infty }\frac{(-1)^{r} (r+1)!\, (2x)^{2\, r+1} }{\left[2\, (r+1)\right]!} .
\end{equation}
The above functions \eqref{GrindEQ__4_}-\eqref{GrindEQ__5_} will be referred as cosine and sine Gaussian, respectively. 
\end{prop}
\begin{proof}
By the use of standard algebraic manipulations, Eq. \eqref{GrindEQ__1_} becomes
\begin{equation*}
e^{-x^{2} } =\frac{1}{1+\hat{c}\, x^{2} } \varphi _{0} = \frac{1}{2} \left(\frac{1}{1-i\, \hat{c}^{\frac{1}{2} } x} +\frac{1}{1+i\, \hat{c}^{\frac{1}{2} } x} \right)\, \varphi _{0}  :=C_g(x)
\end{equation*}
and the associated function $S_{g} (x)$ is 
\begin{equation}\label{GrindEQ__5bis_} 
\begin{split}
S_{g} (x):&=\frac{1}{2\, i} \left[\frac{1}{1-i\, \hat{c}^{\frac{1}{2} } x} -\frac{1}{1+i\, \hat{c}^{\frac{1}{2} } x} \right]\, \varphi _{0} =
\sum _{r=0}^{\infty } (-1)^r  \hat{c}^{r+\frac{1}{2}}x^{2r+1}\varphi_{0}=
\\
& =\sum _{r=0}^{\infty }\frac{(-1)^{r} x^{2\, r+1} }{\Gamma \left(r+\frac{3}{2} \right)} 
 =\frac{1}{\sqrt{\pi } } \sum _{r=0}^{\infty }\frac{(-1)^{r} r!\, (2x)^{2\, r+1} }{\left(2r+1\right)!} =e^{-x^2}Erfi(x) .
\end{split}
\end{equation}
It is worth noting that the Gaussian sine does not represent something new in the scenario of special functions and that it can be indeed interpreted in terms of the Dawson integral $F(x)$, namely \cite{Nijimbere}
\begin{equation}\label{Daw}
\begin{split}
& S_g(x)=\frac{2}{\sqrt{\pi } }\;F(x),\\
& F(x)=e^{-x^2}\int_0^x e^{y^2}dy=\frac{\sqrt{\pi}}{2}e^{-x^2}Erfi(x)=\sum_{n=0}^\infty \frac{(-1)^n2^nx^{2n+1}}{(2n+1)!!}.
\end{split}
\end{equation}
In passing we noted that the "Gaussian trigonometric identity" is
\begin{equation}\label{Gti}
S_g(x)^2+C_g(x)^2=e^{-2x^2}\left(1+Erfi(x)^2 \right). 
\end{equation}
\end{proof}	
In Fig. \ref{fig1} we show $C_g(x)$ vs. $S_g(x)$ resulting in an egg like form, charaterizing the "circle" of the Gaussian trigo and in Fig. \ref{fig2} we also report the behavior of $C_g(x)$ and $S_g(x)$  vs $x$.

\section{Higher Order Gaussian Trigonometric Functions}

After the introductory remarks of the previous section, we introduce higher-order forms of Gauss trigonometric functions. To accomplish this task we take advantage of well known Laplace transform identities \cite{Abramovitz}.

\begin{prop}
 The use of the Laplace transform ($LT$) technique allows to write
\begin{equation}\label{eqs10}
 C_{g} (x)=\int _{0}^{\infty }e^{-s}  \cos \left( \sqrt{\hat{c}} \;x\, s\right)  ds\, \varphi _{0}\; , \;\;
 S_{g} (x)=\int _{0}^{\infty }e^{-s}  \sin \left( \sqrt{\hat{c}} \;x\, s\right)  ds\varphi _{0} .
\end{equation}      
\end{prop}     
\begin{proof}
	The first equation simply follows by the Laplace transform $\frac{k}{k^2+a^2}=\int_0^\infty e^{-kt}\cos(at)dt$
	\begin{equation*}
C_{g} (x)=\dfrac{1}{1+(\sqrt{\hat{c}}\;x)^2}\varphi_{0}=\int_0^\infty e^{-s}\cos\left(\sqrt{\hat{c}}\;x\;s\right) ds\;\varphi_{0}.
	\end{equation*}
The derivation of the second of Eqs. \eqref{eqs10} follows analogous lines.
	\end{proof}
\begin{cor}
The successive derivatives of the Gaussian trigonometric functions can be written in Umbral form as
\begin{equation}
\begin{split}
& C_{g}^{(m)} (x)=\sqrt{\hat{c}^m}\int _{0}^{\infty }e^{-s}  s^{m} \cos \left( \sqrt{\hat{c}} \;x\, s+m\frac{\pi }{2} \right)  ds\, \varphi _{0}\; , \\ 
& S_{g}^{(m)} (x)=\sqrt{\hat{c}^m} \int _{0}^{\infty }e^{-s}  s^{m} \sin \left( \sqrt{\hat{c}} \;x\, s+m\frac{\pi }{2} \right)  ds\, \varphi _{0}\; .
\end{split}  
\end{equation}  
\end{cor}           

\begin{cor}     
The Hermite polynomials are associated with the successive derivatives of the Gaussian $H_{n} (x) e^{-x^{2} } =(-1)^n \partial_x^n e^{-x^2}$, therefore we obtain the identity  
\begin{equation} \label{GrindEQ__8_} 
H_{m} (x)\, \, e^{-x^{2} } =(-1)^{m} \sqrt{\hat{c}^m} \int _{0}^{\infty }e^{-s}  s^{m} \cos \left( \sqrt{\hat{c}} \;x\, s+m\frac{\pi }{2} \right) ds\, \varphi _{0}  \;,
\end{equation} 
which can be exploited for an alternative definition of this family of polynomials.
\end{cor} 

\noindent Let us now go back to the definition of the Gaussian sine which, by using Eq. \eqref{GrindEQ__5bis_}, can also be specified as
\begin{equation} \label{GrindEQ__9_} 
S_{g} (x)=\frac{\hat{c}^{\frac{1}{2} } \, x}{1+\, \hat{c}x^{2} } \, \varphi _{0} \;,
\end{equation} 
and provide the integral in the following example.

\begin{exmp}
	$\forall x \in\mathbb{R}$ 
	\begin{equation}\label{GrindEQ__10_} 
	\int_0^x S_{g} (\xi) \, d\xi =
	 -\dfrac{1}{4\sqrt{\pi}}\sum_{r=1}^\infty \dfrac{(-1)^r (2x)^{2r} (r-1)!}{r\;(2r-1)!}.
	\end{equation}
	indeed
\begin{equation} 
\begin{split}
\int_0^x S_{g} (\xi) \, d\xi &=\frac{\hat{c}^{-\frac{1}{2} } }{2} \ln (1+\hat{c}\, x^{2} )\, \varphi _{0} =\frac{\hat{c}^{-\frac{1}{2} } }{2}(-1)\sum_{r=1}^\infty \dfrac{(-1)^r (\hat{c}x^2)^r}{r}\varphi_{0}= \\
& =-\frac{1}{2}\sum_{r=1}^\infty \dfrac{(-1)^r x^{2r}}{r\;\Gamma\left( r+\frac{1}{2}\right) }=  -\dfrac{1}{4\sqrt{\pi}}\sum_{r=1}^\infty \dfrac{(-1)^r (2x)^{2r} (r-1)!}{r\;(2r-1)!}.
\end{split}
\end{equation} 
\end{exmp}

\begin{rem}
The Gaussian trigonometric functions can be accordingly derived from the real and imaginary parts of the function
\begin{equation}\label{Eg}
E_{g} (x)=\frac{1+i\, \hat{c}^{\frac{1}{2} } \, x}{1+\, \hat{c}x^{2} } \, \varphi _{0} \;.  
\end{equation}    
\end{rem}   
\begin{obs}                             
It is worth stressing that the two Gaussian trigonometric functions are linked by the Kramers-Kronig identity \cite{Jackson} 
\begin{equation}
S_{g} (x)=-\frac{1}{\pi } \mathcal{P}\int _{-\infty }^{\infty }\frac{C_{g} (\xi )}{\xi -x}  d\xi  ,
\end{equation}    
where $\mathcal{P}$ is the Cauchy principal value\footnote{We remind that the Cauchy principal value is " $\mathcal{P}\int _{-\infty }^{\infty }\frac{C_{g} (\xi )}{\xi -x}  d\xi=\lim\limits_{\varepsilon\rightarrow 0^+}\left( \int _{-\infty }^{x-\varepsilon }\frac{C_{g} (\xi )}{\xi -x}  d\xi+\int _{x+\varepsilon}^{\infty }\frac{C_{g} (\xi )}{\xi -x}  d\xi\right) $ ".}. 
\end{obs}             

\begin{figure}[h]
	\centering
	\includegraphics[width=0.5\linewidth]{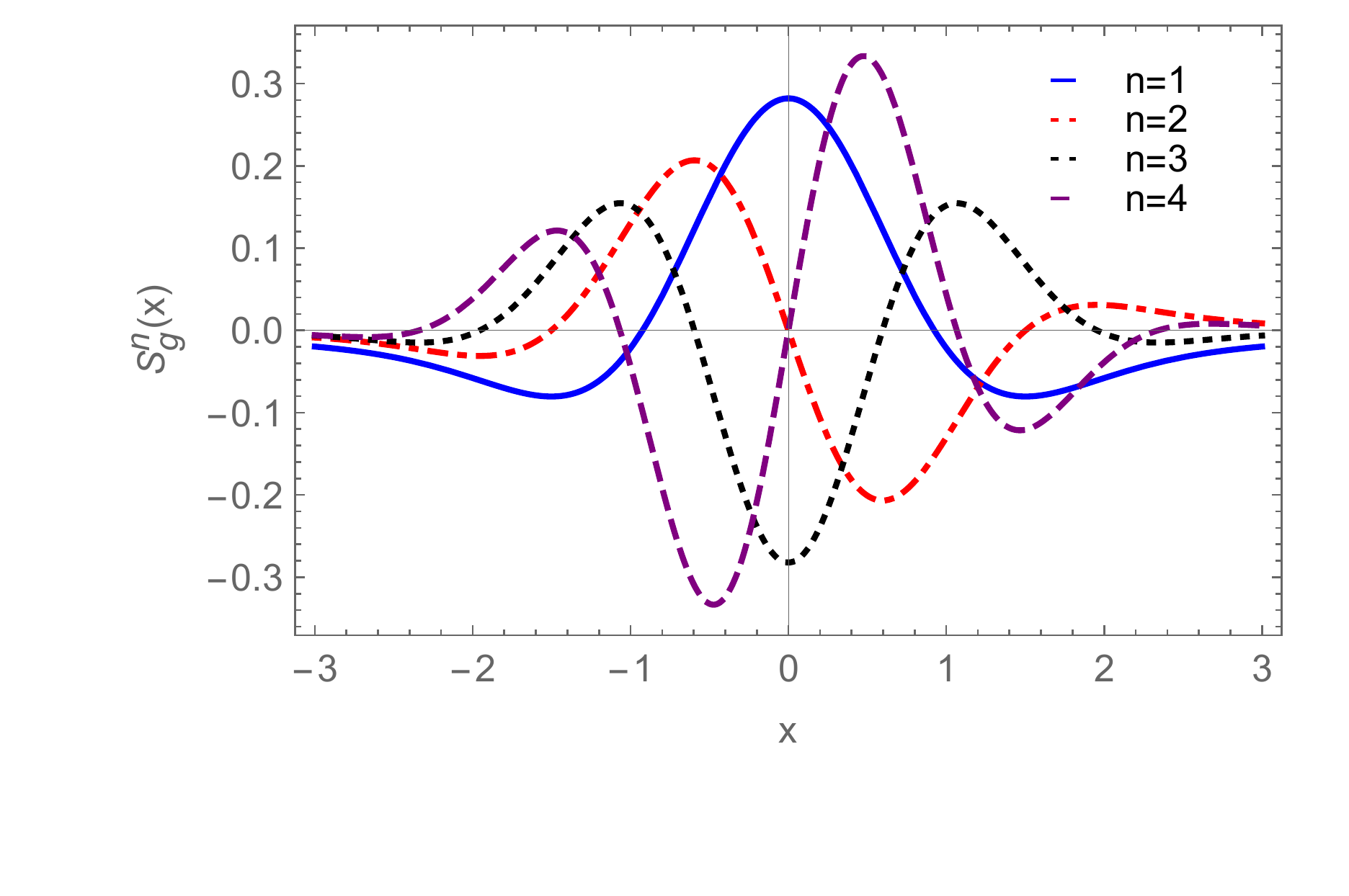}
	\caption{Normalized $S^{(n)}_{g} (x)$ for different values of $n$.}\label{fig3} 
\end{figure}

\begin{prop}                                        
The successive derivative of the complex function $S_{g} (x)$ can be derived from\footnote{This is an interesting result, because the $E_g(x)$, as commented in the final section, is equivalent to the plasma dispersion function.} those of $E_{g} (x)$
\begin{equation} \label{GrindEQ__13_} 
S_g^{(m)}(x)=\frac{2^{m+1} }{\sqrt{\pi } } \sum _{r=\lceil\frac{m-1}{2}\rceil }^{\infty }\frac{(-1)^{r} (r+1)!\, (2\, r+1)!\, (2x)^{2\, r+1-m} }{\left[2\, (r+1)\right]\, !\left(2\, r+1-m\right)\, !}  .
\end{equation} 
\end{prop} 
\begin{proof}
	By exploiting Eq. \eqref{Eg} and the identity based on the Laplace transform of the Lorentz function, we get
	\begin{equation}\label{key}
	E_{g}^{(n)}\! (x)=\partial _{x}^{n} \int _{0}^{\infty }(1+i \hat{c}^{\frac{1}{2} }  x) e^{-s} e^{-s \hat{c} x^{2} } ds \varphi _{0} = Re\left[E_{g}^{(n)} (x)\right]+i Im\left[E_{g}^{(n)} (x)\right]\!.
	\end{equation}
	We study the individual parts by using the two variable Hermite polynomials \cite{Appel,Babusci}
	\begin{equation}\label{key}
	H_{n} (x,\, y)=n!\, \sum _{r=0}^{\lfloor\frac{n}{2} \rfloor}\frac{x^{n-2\, r} y^{r} }{(n-2\, r)!\, r!}
	\end{equation}
	and their properties 
	\begin{align} 
	\partial _{x}^{n} \;e^{a\, x^{2} } =H_{n} (2\, a\, x,\, a)e^{a\, x^{2} },\label{GrindEQ__14_} \\[1.1ex]
	\sum _{k=0}^{\infty }\frac{t^{k} }{k!}  H_{k} (x,\, y)=e^{x\, t+y\, t^{2} }\label{GrindEQ__14b_}  
	\end{align} 	
	so obtaining
	\begin{equation}
	\begin{split}
	Re\left[E_{g}^{(n)} (x)\right]&=C_g^{(n)}(x)=\int _{0}^{\infty} e^{-s}\partial_x^n e^{-s\;\hat{c}\;x^2}ds\varphi_{0}=\\
	& =
	\int _{0}^{\infty} e^{-s}  H_{n} (-2\, s\, \hat{c}\, x,\, -s\, \hat{c})\, e^{-s\, \hat{c}\, x^{2} } ds\varphi _{0}
		\end{split}
	\end{equation}
	and
	\begin{equation} \label{GrindEQ__13b_} 
	\begin{split}
	& Im\left[E_{g}^{(n)} (x)\right]=S_g^{(n)}(x)=\hat{c}^{\frac{1}{2} }\int_0^{\infty} e^{-s} \partial_x^n \left( x\; e^{-s\;\hat{c}\;x^2}\right) ds\varphi_{0}=\\
	& = \hat{c}^{\frac{1}{2} } \sum _{r=0}^{n}\binom{n}{r}\left(\partial _{x}^{n-r} x\right)\int _{0}^{\infty }e^{-s}  H_{r} (-2\, s\, \hat{c}\, x,\, -s\, \hat{c})\, e^{-s\, \hat{c}\, x^{2} } ds\varphi _{0} =\\
	& =\frac{2^{n+1} }{\sqrt{\pi } } \sum _{r=\lceil\frac{n-1}{2}\rceil }^{\infty }\frac{(-1)^{r} (r+1)!\, (2\, r+1)!\, (2x)^{2\, r+1-n} }{\left[2\, (r+1)\right]\, !\left(2\, r+1-n)\right)\, !}  .
	\end{split} 
	\end{equation} 
\end{proof}
The behavior of the successive derivatives of the Gaussian sine function is reported in Fig. \ref{fig3}.

\section{New Forms of Gaussian-like Distribution}

We extend here the methods we just outlined in the previous sections introducing new Gauss like functions and studing the relevant properties.
The use of the methods we have just touched displays a wide flexibility as we are going to show in the examples below.

\begin{exmp}
 We define the function by means of the series
\begin{equation} \label{GrindEQ__15_} 
e\left(-x^{2} |\, n\right):=\sum _{r=0}^{\infty}\frac{(-1)^{r} \Gamma \left(\frac{r}{n} +1\right)\, x^{2\, r} }{r!}, \qquad \forall x\in\mathbb{R}, \forall n\in\mathbb{N}:n>1   
\end{equation} 
and ask whether its integral 
\begin{equation} \label{GrindEQ__16_} 
I_e(n):=\int _{-\infty }^{\infty }e(-x^{2}  |\, n)\, dx 
\end{equation} 
can be calculated in analytical form. This can be checked straightforwardly by the use of the umbral technique discussed so far\footnote{Conventional methods based on the Laplace transform technique can be used too but they are more involute.}. We first write the function \eqref{GrindEQ__15_} in a Gaussian like form through the $\hat{p}$ umbral operator and the following umbral image\footnote{It is easily proved that \begin{equation*}\label{key}
e^{-\hat{p}\, x^{2} } \gamma _{0}=\sum_{r=0}^\infty \dfrac{(-1)^rx^{2r} \hat{p}^r}{r!}\gamma_{0}=_{Eq.  \eqref{uc}}\sum_{r=0}^\infty \dfrac{(-1)^rx^{2r}}{r!} \Gamma\left(\dfrac{r}{n}+1 \right)= e\left(-x^{2} |\, n\right)	.
	\end{equation*}} 
\begin{equation}\label{uc}
 e\left(-x^{2} |\, n\right)=e^{-\hat{p}\, x^{2} } \gamma _{0} , \qquad\qquad \hat{p}^{\alpha } \gamma_{0} =\Gamma \left(\frac{\alpha }{n} +1\right).
\end{equation}
 then we are allowed to evaluate the integral \eqref{GrindEQ__16_} as 
\begin{equation}\label{GrindEQ__18_}
I_e(n)=\int _{-\infty }^{\infty }e^{-\hat{p}\, x^{2} } dx\;\gamma _{0} =
 \sqrt{\frac{\pi }{\hat{p}} } \, \gamma_{0} =\Gamma \left(1-\frac{1}{2\, n} \right)\, \sqrt{\pi } .
\end{equation}
 We note that by taking $n=1$ in the previous equations, namely for $e\left(-x^{2} |\, 1\right)=\frac{1}{1+x^{2}}$, Eq. \eqref{GrindEQ__18_} yields, as it must be, $I_e(1)=\pi $.\\
 
 The result reported in Eq. \eqref{GrindEQ__18_} has been as usual checked numerically. The termwhise integration, we have used as benchmark, becomes difficult. The corresponding series converges very slowly for small values of $n$ and the integration becomes unstable.
We have therefore tested Eq. \eqref{GrindEQ__18_} by the use of the Ramanujian Master Theorem \cite{Berndt,Rama} and by other means commented below.
\end{exmp}

\begin{exmp}
The same procedure yields the following further example 
\begin{equation}\label{GrindEQ__18b_}
\begin{split}
I_{e^{a,b}}(n)&=\int _{-\infty }^{\infty } e(-ax^2+bx \mid n)dx=
\int _{-\infty }^{\infty }e^{-\hat{p}\left( ax^{2}-bx \right) } dx\;\gamma _{0} =\\
& =
\sqrt{\frac{\pi }{a\hat{p}} }e^{\hat{p}\frac{b2}{4a}} \gamma_{0}  = \sqrt{\frac{\pi }{a} } \sum_{r=0}^\infty \dfrac{\left( \frac{b^2}{4a}\right) ^r}{r!}\hat{p}^{r-\frac{1}{2}}\gamma_{0}=
 \sqrt{\frac{\pi }{a} } \;e_{-\frac{1}{2}}\left( \dfrac{b^2}{4a}\mid n\right) ,\\
e_\nu(x\mid n)&= \sum_{r=0}^\infty \dfrac{x^r}{r!}\Gamma\left(\dfrac{r+\nu}{n}+1 \right) .
\end{split}
\end{equation}
It is finally interesting to note that the use of the generating function

\begin{equation}\label{key}
\begin{split}
& \sum_{s=0}^\infty \dfrac{t^s}{s!}H_s(x,y\mid n)= e(xt+yt^2 \mid n), \\
& H_m(x,y\mid n) =m! \sum_{r=0}^{\lfloor \frac{m}{2}\rfloor} \dfrac{x^{m-2r}y^r}{(m-2r)!r!}\Gamma\left( \dfrac{m-r}{n}+1\right) 
\end{split}
\end{equation}
allows the introduction of a further family of two-variable Hermite like polynomials, whose properties will be explored elsewhere.
\end{exmp}

\begin{exmp}
It is evident that the function \eqref{GrindEQ__15_}, which with increasing $n$ becomes closer and closer to a Gaussian, can be exploited to model distributions with a tail longer than an ordinary Gauss function. By setting indeed
\begin{equation}
F(x;\, \sigma |n):=\frac{1}{\sqrt{2\, \pi } \, \Gamma \left(1-\frac{1}{2\, n} \right) \sigma } e\left(-\frac{x^{2} }{2\, \sigma ^{2} } |n\right) ,
\end{equation}
we can evaluate the moments associated with its distribution by using the generating function method \eqref{GrindEQ__14b_}. To this aim we note that, $\forall d\in\mathbb{N}, \;\forall \sigma\in\mathbb{R}:\sigma\neq 0$,
\begin{equation} \label{GrindEQ__19_} 
\langle F(x;\, \sigma |n) \rangle_{(m,\;d)} :=M_{(m,\;d)} =\int _{-\infty }^{\infty }(x+d)^{m}  F(x,\, \sigma |n)\, dx \;, \qquad \forall m\in\mathbb{Z},
\end{equation} 
can be calculated by noting that
\begin{equation} \label{GrindEQ__20_} 
\begin{split}
\sum _{m=0}^{\infty }\frac{t^{m} }{m!}  M_{(m,\;d)} &=
\int _{-\infty }^{\infty }\sum _{m=0}^{\infty }\frac{t^{m} (x+d)^m}{m!}F(x,\, \sigma |n)\, dx
=\\
& =\frac{e^{t\;d} }{\sqrt{2\, \pi } \, \Gamma \left(1-\frac{1}{2\, n} \right)\sigma } \int _{-\infty }^{\infty }e^{t\, x}  e^{-\frac{\hat{p}\, x^{2} }{2\, \sigma ^{2} } } \, dx\, \gamma_{0}  .
\end{split}
\end{equation} 
The Gaussian integral on the r.h.s. of Eq. \eqref{GrindEQ__20_} yields 
\begin{equation} \label{GrindEQ__21_} 
\sum _{m=0}^{\infty }\frac{t^{m} }{m!}  M_{(m,\;d)} =
 \frac{e^{td}\;e^{\frac{t^{2} \sigma ^{2} }{2\, \hat{p}} } }{\Gamma \left(1-\frac{1}{2\, n} \right)} \hat{p}^{-\frac{1}{2} } \gamma _{0}
\end{equation} 
and the use of the Hermite generating function \eqref{GrindEQ__14_} 
finally provides the result 
\begin{equation} \label{GrindEQ__23_} 
\begin{split}
M_{(m,\;d)} &=\, \frac{\hat{p}^{-\frac{1}{2} } }{\Gamma \left(1-\frac{1}{2\, n} \right)} H_{m} \left( d,\, \frac{\sigma ^{2} }{2} \hat{p}^{-1} \right)  \gamma_{0}  =\\
& =\frac{m!}{\Gamma \left(1-\frac{1}{2\, n} \right)} \sum _{r=0}^{\lfloor\frac{m}{2} \rfloor}\frac{ d^{m-2\, r} \sigma ^{2\, r} }{2^{r} r!\, (m-2\, r)!}  \Gamma \left(1-\dfrac{1}{2n}-\frac{r }{n} \right).
\end{split}
\end{equation} 
The higher order moments for $d=0$ can accordingly be evaluated as 
\begin{equation}
\langle F(x; \sigma |n) \rangle_m :=M_m =\frac{1}{\Gamma \left(1-\frac{1}{2 n} \right)} \frac{m!}{\Gamma\left( \frac{m}{2}+1\right) } \left(\frac{\sigma ^{2} }{2} \right)^{\!\frac{m}{2} } \Gamma \left(1-\dfrac{1}{2n}-\frac{m}{2 n} \right)\!.
\end{equation}
The quasi Gaussian distribution have only a finite number of non diverging higher order moments, compatible with the condition $\frac{m+1}{2\, n} <1$. These distributions can be exploited to interpolate between Gaussian and Cauchy-Lorentz distributions.
\end{exmp} 

\begin{exmp}
 It is worth noting that the use of the integral representation of the Gamma function allows to cast Eq. \eqref{GrindEQ__15_} in the form
\begin{equation}\label{GrindEQ__25_}
e\left(-x^{2} |\, n\right)=\sum _{r=0}^{\infty}\frac{(-1)^{r} x^{2\, r} }{r!}  \int _{0}^{\infty }e^{-s} s^{\frac{r}{n} }  ds=\int _{0}^{\infty }e^{-s}  e^{-x^{2} s^{\alpha } } ds,\qquad \alpha =\frac{1}{n} .
\end{equation}
By expoliting the property \eqref{GrindEQ__14_}, the successive derivatives of $e\left(-x^{2} |\, n\right)$ can be then written as
\begin{equation}\label{alp}
e^{\left(m\right)} \left(-x^{2} |\, n\right)=\int _{0}^{\infty }e^{-s} H_{m} (-2\, x\, s^{\alpha } ,\,- s^{\alpha } )\,  e^{-x^{2} s^{\alpha } } ds,\qquad \alpha =\frac{1}{n} 
\end{equation}
and can be exploited to establish families of functions providing a smooth transition, with increasing $n$ , to the ordinary Hermite Gauss functions.\\

\noindent The integral transform \eqref{GrindEQ__25_} is an alternative to the umbral notation developed so far and show noticeable features of interest, which will be touched below and will be discussed more carefully in a forthcoming investigation. 
\end{exmp}

It is also worth noting that the use of the integral representation in Eq. \eqref{GrindEQ__25_} yields almost straightforwardly the evaluation of its infinite integral as shown below.

\begin{exmp}
 We can accordingly write
\begin{equation}\label{key}
I(\alpha):=\int_{-\infty }^\infty \int_0^\infty e^{-s-x^2s^\alpha}ds\;dx.
\end{equation}
Interchanging the integrals and by using the ordinary Gaussian 
and the properties of the Gamma function \cite{Abramovitz}, when $Re(\alpha)<2$, we end up with
\begin{equation}
\begin{split}
I(\alpha)&=\int_0^\infty e^{-s}\int_{-\infty}^\infty e^{-x^2s^\alpha} dx\;ds=
\int_0^\infty e^{-s} \sqrt{\dfrac{\pi}{s^\alpha}}ds=\\
& =
\sqrt{\pi}\int_0^\infty e^{-s} s^{-\frac{\alpha}{2}}ds=\sqrt{\pi}\; \Gamma\left(1-\dfrac{\alpha}{2} \right) 
\end{split}
\end{equation}
and the same procedure can be exploited for the derivation of higher order moments.
\end{exmp}

\section{Applications and Final Comments}

It is evident that the discussion on the properties of Gauss trigonometric functions can be extended to the case of the distribution in Eq. \eqref{GrindEQ__15_}. We will not discuss these problems any more, leaving the subject for a  forthcoming investigation. We will  devote this concluding section to the relevance of the previous discussion to the link with known families of special functions.\\

We have already mentioned that the Gaussian sine is linked to the Dawson function, through Eq. \eqref{Daw}, keeping advantage from the relevant expression in terms of confluent hypergeometric series ${}_1F_1(.;.;.)$ \cite{Nijimbere2}
\begin{equation}\label{SgF}
S_g(x)=\dfrac{2}{\sqrt{\pi}}\;e^{-x^2}\int_0^x e^{y^2}dy=\dfrac{2}{\sqrt{\pi}}\;x\;e^{-x^2}{}_1F_1\left(\dfrac{1}{2};\dfrac{3}{2};x^2 \right) .
\end{equation}
This identity offers a further direction along which we can extend the umbral point of view.

\begin{exmp}
We can consider indeed the following umbral interpretation of the confluent hypergeometric function ${}_1F_1(a;c;x)$ through the Pochhammer symbol $(y)_r=y(y+1)\cdots(y+r-1)$
\begin{equation}\label{Hyp}
\begin{split}
& {}_1F_1(a;c;x)=\sum_{n=0}^\infty \dfrac{(a)_n\; x^n}{(c)_n \;n!}=\sum_{n=0}^\infty \hat{\kappa}^n\dfrac{x^n}{n!}\phi_0=e^{\hat{\kappa}x}\phi_0,\\
& \hat{\kappa}^n\phi_0:=\hat{\kappa}_{a,c}^n\;\phi_0 =\dfrac{(a)_n}{(c)_n}. 
\end{split}
\end{equation}
For umbral and Pochhammer properties, we note that
\begin{equation}\label{key}
\hat{\kappa}^n\hat{\kappa}^m\phi_0=\hat{\kappa}^{n+m}\phi_0=\dfrac{(a)_{n+m}}{(c)_{n+m}}=\dfrac{(a)_n}{(c)_n}\dfrac{(a+n)_m}{(c+n)_m}.
\end{equation}
By using above relation, we can easily calculate high order of derivative of hypergeometric function ${}_1F_1(a;c;x)$
\begin{equation}
\begin{split}
\partial_x^s\; {}_1F_1(a;c;x)&=\hat{\kappa}^s e^{\hat{\kappa}x}\phi_0=\hat{\kappa}^s\sum_{n=0}^\infty \hat{\kappa}^n \dfrac{x^n}{n!}\phi_0=\sum_{n=0}^\infty\dfrac{(a)_s}{(c)_s}\dfrac{(a+s)_n}{(c+s)_n}\dfrac{x^n}{n!}=\\
& =\dfrac{(a)_s}{(c)_s}\;{}_1F_1(a+s;c+s;x).
\end{split}
\end{equation}
According to Eq. \eqref{Hyp}, if we get $a=\frac{1}{2}$ and $c=\frac{3}{2}$, Eq. \eqref{SgF} becomes 
\begin{equation}\label{Sgk}
S_g(x)=\dfrac{2}{\sqrt{\pi}}\;x\;e^{-(1-\hat{\kappa})x^2}\phi_0.
\end{equation}
Namely, a fairly straightforward form which can be usefully applied to perform specific calculations, involving this family of functions.
\end{exmp}
Regarding integrals involving the the gaussian sine we provide some examples, as shown below.

\begin{exmp}
From Eq. \eqref{Sgk} we get, if $\mid\alpha\mid<1$, 
\begin{equation}\label{key}
\begin{split}
1) \;& \int_{0}^\infty S_g(x,\alpha)dx=\dfrac{2}{\sqrt{\pi}}\int_0^\infty x\;e^{-(1-\alpha\hat{\kappa})x^2}dx\;\phi_0=
\dfrac{1}{\sqrt{\pi}}\dfrac{1}{1-\alpha\hat{\kappa}}\phi_0=\\
& =\dfrac{1}{\sqrt{\pi}}\sum_{r=0}^\infty \alpha^r\hat{\kappa}^r\phi_0 =\dfrac{1}{\sqrt{\pi}}\sum_{r=0}^\infty \dfrac{\left(\frac{1}{2} \right)_r }{\left(\frac{3}{2} \right)_r }\alpha^r= \dfrac{1}{\sqrt{\pi}}\;{}_2F_1\left(\dfrac{1}{2},1;\dfrac{3}{2};\alpha \right).
\end{split}
\end{equation}	

\begin{equation}\label{key}
\begin{split}
2) \;&\int_{-\infty}^\infty \dfrac{S_g(x)}{x}dx=\dfrac{2}{\sqrt{\pi}}\int_{-\infty}^\infty e^{-(1-\hat{\kappa})x^2}dx\;\phi_0=
\dfrac{2}{\sqrt{1-\hat{\kappa}}}\phi_0=\\
& =\!
2\sum_{r=0}^\infty \binom{r\!-\!\frac{1}{2}}{r}\hat{\kappa}^r\phi_0=2\sum_{r=0}^\infty \binom{r-\frac{1}{2}}{r} \dfrac{\left(\frac{1}{2} \right)_r }{\left(\frac{3}{2} \right)_r } = 2\cdot {}_2F_1\left(\dfrac{1}{2}, \dfrac{1}{2};\dfrac{3}{2};1\right)\! =\!\pi.
\end{split}
\end{equation}	
\end{exmp}
\begin{exmp}
The Fresnel integrals can be written in terms of hypergeometric functions ${}_1F_2(a;b,c;x)$ as \cite{Nijimbere2} 
\begin{equation}\label{key}
\begin{split}
& C(x)=\int_0^x \cos\left( \dfrac{\pi}{2} \eta^2\right) d\eta=x \;{}_1F_2\left(\dfrac{1}{4};\dfrac{1}{2},\dfrac{5}{4};-\left( \dfrac{\pi}{4}x^2\right) ^2 \right) ,\\
& S(x)=\int_0^x \sin\left( \dfrac{\pi}{2} \eta^2\right)d\eta=\dfrac{\pi}{2}\; \dfrac{x^3}{3} \;{}_1F_2\left(\dfrac{3}{4};\dfrac{3}{2},\dfrac{7}{4};-\left( \dfrac{\pi}{4}x^2\right) ^2 \right) ,\\
& {}_1F_2(a;b,c;x)=\sum_{n=0}^\infty \dfrac{(a)_n}{(b)_n (c)_n}\dfrac{x^n}{n!}.
\end{split}
\end{equation}
According to the previous umbral discussion, the relevant images are Gauss functions and therefore, we find
\begin{equation}\label{key}
\begin{split}
& C(x)=x\;e^{-\left( \frac{\pi}{4}x^2\right)^2 \hat{\chi}_c}\;\phi_{0,c}, \qquad\qquad \hat{\chi}_c^n:={}_\frac{1}{4}\hat{\chi}_{\left( \frac{1}{2},\frac{5}{4}\right) }\phi_{0,c}=\dfrac{\left( \frac{1}{4}\right)_n }{\left(\frac{1}{2} \right)_n\left(\frac{5}{4} \right)_n  }\\
& S(x)=\pi\; \dfrac{x^3}{3}\;e^{-\left( \frac{\pi}{4}x^2\right)^2 \hat{\chi}_s}\;\phi_{0,s}, \qquad\qquad \hat{\chi}_s^n:={}_\frac{3}{4}\hat{\chi}_{\left( \frac{3}{2},\frac{7}{4}\right) }\phi_{0,s}=\dfrac{\left( \frac{3}{4}\right)_n }{\left(\frac{3}{2} \right)_n\left(\frac{7}{4} \right)_n  }
\end{split}
\end{equation}
\end{exmp}

\begin{exmp}
It is evident that the previous results can be exploited to simplify computations involving this family of functions.
For example, if we are interested in computing improper integrals involving the Fresnel functions, as e.g.
\begin{equation}\label{key}
\int_0^\infty \dfrac{S(\xi)}{\xi^3}d\xi=\dfrac{\pi}{3}\int_0^\infty e^{-\left( \frac{\pi}{4}\xi^2\right)^2 \hat{\chi}_c}\;d\xi\;\phi_{0,s},
\end{equation}
we can easily achieve our result, by recalling that 
\begin{equation}\label{key}
\int_0^\infty e^{-ax^4}dx=\dfrac{1}{4}\Gamma\left(\dfrac{1}{4}\right) a^{-\frac{1}{4}}, \qquad Re(a)>0
\end{equation}
and, by replacing the parameter $"a"$ with $\left( \frac{\pi}{4}\xi^2\right)^2 \hat{\chi}_c$ we eventually end up with 
\begin{equation}
\begin{split}
\dfrac{\pi}{3}\int_0^\infty e^{-\left( \frac{\pi}{4}\xi^2\right)^2 \hat{\chi}_c}\;d\xi\;\phi_{0,c}&= 
\dfrac{\sqrt{\pi}}{6} \Gamma\left(\dfrac{1}{4}\right)\hat{\chi}^{-\frac{1}{4}}\phi_{0,s}=\\
& =
\dfrac{\sqrt{\pi}}{6} \Gamma\left(\dfrac{1}{4}\right)\dfrac{\left( \frac{3}{4}\right)_{-\frac{1}{4}} }{\left(\frac{3}{2} \right)_{-\frac{1}{4}}\left(\frac{7}{4} \right)_{-\frac{1}{4}}  }.
\end{split}
\end{equation}
\end{exmp}

The applicative framework of the results we have just obtained is indeed interesting. We like to mention the Fried-Conte dispersion function $Z(x)$, often used in Plasma Physics \cite{Fried} which, within the present context, are just written in terms of the complex function defined in Eq. \eqref{Eg} 
\begin{equation}\label{k.ey}
Z(x)=i\sqrt{\pi}E_g(x).
\end{equation}

Before closing the paper, we consider interesting to add a comment, regarding the umbral definition of the so-called Lévy stable distributions \cite{Kahane,Anderssen,PG}, describing non standard statistical effects in different phenomenological environments.
Regarding the Umbral form of this family of functions, we proceed as outlined below.

 \begin{exmp}
 L\'evy stable distribution $g_{\alpha } (x)$, in the present notation, can be defined as \cite{SLicciardi}
\begin{equation} \label{GrindEQ__31_} 
g_{\alpha } (x)=-\frac{1}{\pi \, x} e^{-\left(\frac{\hat{f}}{x} \right)^{\alpha } } \varepsilon_{0} , \;\; \forall \alpha\in\mathbb{R}:0<\alpha <1, \;\;
 \hat{f}^{\beta } \varepsilon_{0} =\Gamma (\beta \, +1)\sin (\pi \, \beta )
\end{equation} 
and their explicit expression in terms of infinite series writes indeed 
\begin{equation}
g_{\alpha } (x)=-\frac{1}{\pi } \sum _{r=0}^{\infty }\frac{(-1)^{r} x^{-\alpha \, r-1} }{r!}  \Gamma (\alpha \, r+1)\, \sin \left(\pi \, \alpha \, r\right) .
\end{equation}                        
Albeit of limited usefulness for accurate numerical computations, as underscored in \cite{PG}, we have used this form because of
its straightforward umbral version. By applying the same procedure leading to Eq. \eqref{GrindEQ__25_} we find the following integral representation 
\begin{equation} \label{GrindEQ__33_} 
g_{\alpha } (x)=-\frac{Im}{\pi } \int _{0}^{\infty }e^{-s\, x-(-1)^{\alpha } s^{\alpha } }  ds .
\end{equation} 
These distributions have the remarkable properties that all the moments $\;\left\langle x^{\mu } \right\rangle =\int _{0}^{\infty }x^{\mu }  g_{\alpha } (x)\, dx\;$ are not defined for $\mu >\alpha $. By a direct application of our method it is indeed easy to check that 
\begin{equation}
\begin{split}
 \left\langle x^{\mu } \right\rangle &=-\frac{1}{\pi } \int _{0}^{\infty }x^{\mu -1}  e^{-\left(\frac{\hat{f}}{x} \right)^{\alpha } } \, dx\, \varepsilon_{0} =
 -\frac{1}{\pi \, \alpha } \Gamma \left(-\frac{\mu }{\alpha } \right)\, \hat{f}^\mu \varepsilon _{0}  =\\
 & =-\frac{1}{\pi\, \alpha } \Gamma \left(-\frac{\mu }{\alpha } \right)\, \Gamma \left(\mu+1 \right)\, \sin \left(\pi \mu \right)=\\
 & =\frac{\Gamma(\mu)\sin(\pi\mu)}{\sin\left(\pi\frac{\mu}{\alpha}\right)\Gamma\left( \frac{\mu}{\alpha}\right)  }, \qquad 0 <\mu <\alpha <1.
\end{split}  
\end{equation}   
 \end{exmp}

\begin{exmp}      
The further important property \cite{PG}
\begin{equation} \label{GrindEQ__35_} 
\int _{0}^{\infty }e^{-p\, x} g_{\alpha }  (x)\, dx=e^{-p^{\alpha } }  , \qquad p>0, \;\; 0<\alpha<1,
\end{equation} 
namely the fact that the Laplace transform of L\`{e}vy stable distribution is the stretched exponential \cite{Kohl} , is easily derived from Eq. \eqref{GrindEQ__33_}, as shown below 
\begin{equation} \label{GrindEQ__36_} 
\begin{split}
\int _{0}^{\infty }e^{-p\, x} g_{\alpha }  (x)\, dx&=-\frac{Im}{\pi } \int _{0}^{\infty }e^{-p\, x}  \left[\int _{0}^{\infty }e^{-s\, x-(-1)^{\alpha } s^{\alpha } }  ds\right]\, dx=\\
& =\frac{Im}{\pi } \int _{0}^{\infty }\frac{e^{-(-1)^{\alpha } s^{\alpha } } }{s+p}  ds=e^{-p^{\alpha } } ,\qquad p>0. 
\end{split}
\end{equation} 
Furthermore, the Laplace transform of the following modified form of $g_{\alpha } (x)$
\begin{equation} \label{GrindEQ__37_} 
g_{\alpha ,\, \nu } (x)=-\frac{Im}{\pi } \int _{0}^{\infty }s^{\nu } e^{-s\, x-(-1)^{\alpha } s^{\alpha } }  ds 
\end{equation} 
yields the Weibull distribution \cite{Tsallis,Tsallis2}
\begin{equation} \label{GrindEQ__38_} 
\int _{0}^{\infty }e^{-p\, x} g_{\alpha ,\alpha -1}  (x)\, dx=p^{\alpha -1} e^{-p^{\alpha } }  .
\end{equation} 
 \end{exmp}

\subsection{\textbf{Final Comments}}
In this article we have obtained a number of results which open interesting perspectives. Further developments will be discussed in a forthcoming paper where we will complete the scenario of the possibilities offered by the extension of the umbral formalism to the hypergeometric functions. In particular we will see how these methods can be extended to the study of Fresnell functions and how the method can be exploited to treat problems in classical Optics.\\

\textbf{Acknowledgements}\\

\noindent The work of Dr. S. Licciardi was supported by an Enea Research Center individual fellowship. Furthermore, the paper was carried out under the auspices of INDAM’s GNFM.\\

%

{}

\end{document}